\documentclass[11pt]{pjm}
\usepackage{amssymb, amsthm, amsmath}

\def\CM{Cohen--Macaulay }
\def\deg{\mathop{\mathrm{deg}}}

\def\div{\mathop{\mathrm{div}}}

\def\Proj{\mathop{\mathrm{Proj}}}
\def\Spec{\mathop{\mathrm{Spec}}}

\def\zero{^{\circ}}

\def\Im{\mathop{\mathrm{Im}}}

\title{Veronese subrings and tight closure}

\author{Anurag K. Singh}

\address{Department of Mathematics, University of Michigan,  
Ann Arbor, MI 48109 \\
{\em Current address:} Department of Mathematics, 
University of Illinois, \\
1409 W. Green Street, Urbana, IL 61801}

\email{\tt singh6@math.uiuc.edu}
\urladdr{\tt http://www.math.uiuc.edu$\sim$singh6}

\newtheorem{thm}{Theorem}[section]
\newtheorem{prop}[thm]{Proposition}
\newtheorem{lem}[thm]{Lemma}
\newtheorem{cor}[thm]{Corollary}

\theoremstyle{definition}
\newtheorem{defn}[thm]{Definition}
\newtheorem{ex}[thm]{Example}

\theoremstyle{remark}
\newtheorem{rem}[thm]{Remark}

\date{23 June, 1998, Revised: 19 October, 1998}

\begin{document}

\begin{abstract} 
We determine when an $\mathbb N$--graded ring has Veronese subrings which are
F--rational or F--regular. The results obtained here give a better
understanding of these properties, and include various techniques of
constructing F--rational rings which are not F--regular 
\end{abstract}

\maketitle

\section{Introduction}

Throughout this paper, all rings are commutative, Noetherian, and have an 
identity element. By a {\it graded ring}, we mean a ring $R=\oplus_{n \ge
0}R_n$, which is finitely generated over a field $R_0=K$.

The theory of {\it tight closure}  was developed by Melvin Hochster and Craig
Huneke in \cite{HHjams}, and has yielded many elegant and powerful results in
commutative algebra and related fields. The theory draws attention to rings
which have the property that all their ideals are tightly closed, called {\it
weakly F--regular} rings, and rings with the weaker property that parameter
ideals are tightly closed, called {\it F--rational} rings. The term {\it
F--regular}\/ is reserved for rings with the property that all their
localizations are weakly F--regular. (The recent work of Lyubeznik and Smith
shows that for graded rings the properties of weakly F--regularity and
F--regularity are equivalent, see \cite{LS}.) These properties turn out to be
of significant importance, for instance the Hochster--Roberts theorem of
invariant theory that direct summands of polynomial rings are Cohen--Macaulay
(\cite{HRinv}), can actually be proved for the much larger class of F--regular
rings.

While the property of F--rationality provides an algebraic analogue of the
notion of rational singularities, F--regularity, in general, is not so well
understood geometrically. One approach is to study the variety $X = \Proj R$
for a graded F--regular ring $R$. The Veronese subrings of $R$ are also
homogeneous coordinate rings for $X$, and so it is interesting to determine
when graded rings have F--rational or F--regular Veronese subrings. The
question regarding F--rational Veronese subrings is easily answered: let
$(R,m,K)$ be a \CM graded domain of dimension $d$, with an isolated singularity
at $m$. We show that there exists a positive integer $n$ such that the Veronese
subring $R^{(n)}$ is F--rational if and only if $[H_m^d(R)]_0=0$. With regard
to F--regular Veronese subrings, we show that if $R$ is a normal ring generated
by degree one elements over a field, then either $R$ is F--regular, or else no
Veronese subring of $R$ is F--regular. This leads us to the question: if
$(R,m,K)$ is a normal graded domain, generated by degree one elements, with an
isolated singularity at $m$, then under what conditions is $R$ F--regular? It
is easily seen that F--regularity forces the {\it $a$--invariant}, $a(R)$, to
be negative. For rings of dimension two (although not in higher dimensions)
this is also a sufficient condition for F--regularity. We construct rings $R$
of dimension $d \ge 3$ with $a(R) = 2-d$ which are not F--regular, while if
$a(R) < 2-d$, Smith has pointed out that $\Proj R$ is a variety of minimal
degree, and $R$ is indeed F--regular. We also construct a rich family of
F--rational rings of characteristic zero, with isolated singularities, which
have no F--regular Veronese subrings.

We would like to point out that although tight closure is primarily a
characteristic $p$ notion, it has strong connections with the study of
singularities of algebraic varieties over fields of characteristic zero.
Specifically, let $R$ be a ring essentially of finite type over a field of
characteristic zero. Then $R$ has rational singularities if and only if it is
of F--rational type, see \cite{Hara, Smratsing}. In the $\mathbb Q$--Gorenstein
case, we have some even more remarkable connections: F--regular type is
equivalent to having log--terminal singularities and F--pure type implies (and
is conjectured to be equivalent to) log--canonical singularities, see
\cite{Smvanish, Walog}.

\section{Preliminaries}

Let $R$ be a Noetherian ring of characteristic $p > 0$. We shall always use
the letter $e$ to denote a variable nonnegative integer, and $q$ to denote the
$e\,$th power of $p$, i.e., $q=p^e$. We shall denote by $F$, the Frobenius 
endomorphism of $R$, and by $F^e$, its $e\,$th iteration, i.e., $F^e(r)=r^q$. 
For an ideal $I=(x_1, \dots, x_n) \subseteq R$, we let $I^{[q]}=(x_1^q, \dots,
x_n^q)$. Note that $F^e(I)R= I^{[q]}$, where $q=p^e$, as always. Let $S$
denote the ring $R$ viewed as an $R$--algebra via $F^e$. Then $S\otimes_R\_$
is a covariant functor from $R$--modules to $S$--modules, and so is a
covariant functor from $R$--modules to $R$--modules! If we consider a map of
free modules $R^n \to R^m$ given by the matrix $(r_{ij})$, applying $F^e$ we
get a map $R^n \to R^m$ given by the matrix $(r_{ij}^q)$. For an $R$--module
$M$, note that the $R$--module structure on $F^e(M)$ is $r'(r\otimes m) = r'r
\otimes m$, and $r' \otimes rm = r'r^q \otimes m$. For $R$--modules $N
\subseteq M$, we use $N_M^{[q]}$ to denote $\Im (F^e(N) \to F^e(M))$.

For a reduced ring $R$ of characteristic $p > 0$, $R^{1/q}$ shall denote the
ring obtained by adjoining all $q\,$th roots of elements of $R$. The ring $R$
is said to be {\it F--finite}\/ if $R^{1/p}$ is module--finite over $R$. Note
that a finitely generated algebra $R$ over a field $K$ is F--finite if and only
if $K^{1/p}$ is a finite field extension of $K$.

We shall denote by $R\zero$ the complement of the union of the minimal primes
of $R$. We say $I=(x_1,\dots,x_n) \subseteq R$ is a {\it parameter ideal}\/ if
the images of $x_1,\dots,x_n$ form part of a system of parameters in the local
ring $R_P$, for every prime ideal $P$ containing $I$.

\begin{defn} 
Let $R$ be a ring of characteristic $p$, and $I$ an ideal of $R$. An element
$x$ of $R$, is said to be in $I^F$, the {\it Frobenius closure }\/ of $I$, if
there exists some $q=p^e$ such that $x^q \in I^{[q]}$.

For $R$--modules $N \subseteq M$ and $u \in M$, we say that $u \in N_M^*$, the
{\it tight closure }\/ of $N$ in $M$, if there exists $c \in R\zero$ such that
$cu^q \in N_M^{[q]}$ for all $q=p^e \gg 0$. It is worth recording this when
$M=R$, and $N=I$ is an ideal of $R$. An element $x$ of $R$ is said to be in
$I^*$ if there exists $c \in R\zero$ such that $cx^q \in I^{[q]}$ for all
$q=p^e \gg 0$. If $I^*=I$ we say that the ideal $I$ is {\it tightly closed}. 

A ring $R$ is said to be {\it F--pure }\/ if for all $R$--modules $M$, the 
Frobenius homomorphism $F: M \to F(M)$ is injective. A ring $R$ is {\it weakly
F--regular}\/ if every ideal of $R$ is tightly  closed, and is {\it
F--regular}\/ if every localization is weakly F--regular. Lastly, $R$ is said
to be {\it F--rational}\/ if every parameter ideal of $R$  is tightly closed. 
\end{defn}

It is easily verified that $I \subseteq I^F \subseteq I^*$. Furthermore, $I^*$
is always contained in the integral closure of $I$, and is frequently much
smaller. A weakly F--regular ring is F--rational as well as F--pure. We next
record some useful results.

\begin{thm} 
\item $(1)$\quad Regular rings are F--regular. A ring which is a direct summand 
of an F--regular ring is itself F--regular. 

\item $(2)$\quad An F--rational ring $R$ is normal. If, in addition, $R$ is 
the homomorphic image of a \CM ring, then it is Cohen--Macaulay.

\item $(3)$\quad An F--rational Gorenstein ring is F--regular. 

\item $(4)$\quad Let $(R,m)$ be a reduced excellent local ring of dimension $d$ 
and characteristic $p > 0$. If $c\in R\zero$ is an element such that $R_c$ is
F--rational, then there exists a positive integer $N$ such that 
$c^N(0_{H_m^d(R)}^*)=0$.

\item $(5)$\quad Let $R$ be a graded ring. Then $R$ is weakly F--regular
if and only if it is F--regular.
\label{longlist} 
\end{thm}

\begin{proof} 
For assertions (1)---(3), see \cite[Theorem 4.2]{HHbasec}. Part (4) is a result 
of Velez, \cite{Velez}, and (5) is \cite[Corollary 4.4]{LS}. 
\end{proof}

\begin{rem} 
The equivalence of weak F--regularity and F--regularity, in general, is a 
formidable open question. However in the light of Theorem \ref{longlist} (5)
above, we frequently have no reason to distinguish between these notions. 
\end{rem}

By a graded ring $(R,m,K)$, we shall always mean a ring $R=\oplus_{n \ge 0}R_n$
finitely generated over a field $R_0=K$. We shall denote by $m=R_+$, the
homogeneous maximal ideal of $R$. The {\it punctured spectrum}\/ of $R$ refers
to the set $\Spec R - \{m\}$. By a system of parameters for $R$, we shall mean
a sequence of homogeneous elements of $R$ whose images form a system of
parameters for $R_m$. In specific examples involving homomorphic images of
polynomial rings, lower case letters shall denote the images of the
corresponding variables, the variables being denoted by upper case letters.

For conventions regarding graded modules and homomorphisms, we follow
\cite{GW}. For a graded $R$--module $M$, we shall denote by $[M]_i$, the
$i$-th graded piece of $M$. 

\begin{defn}
Let $R=\oplus_{i\ge 0}R_i$ be a graded ring, and $n$ be a
positive integer. We shall denote by $R^{(n)}$, the {\it Veronese subring}\/
of $R$ spanned by all elements of $R$ which have degree a multiple of $n$,
i.e., $R^{(n)}=\oplus_{i\ge 0}R_{in}$. 
\end{defn}

Note that the ring $R^{(n)}$ is a direct summand of $R$ as an $R^{(n)}$--module
and that $R$ is integral over $R^{(n)}$. Hence whenever $R$ is \CM or normal,
so is $R^{(n)}$. We record the following result, see \cite[Lemme 2.1.6]{EGA} or
\cite[page 282]{Mumford} for a proof.

\begin{lem} 
Let $R$ be a graded ring. Then there exists a positive integer
$n$ such that the Veronese subring $R^{(n)}$ is generated over $K$ by forms
of equal degree. 
\label{existver} 
\end{lem}

Recall that the highest local cohomology module $H_m^d(R)$ of $R$, where $\dim
R =d$, may be identified with $ \varinjlim R/(x_1^t,\dots, x_d^t) $ where $x_1,
\dots, x_d$ is a system of parameters for $R$ and the maps are induced by
multiplication by $x_1 \dotsm x_d$. If $R$ is Cohen--Macaulay, these maps are
injective. The $R$--module $H_m^d(R)$ carries a natural graded structure,
namely $\deg [r+(x_1^t,\dots, x_d^t)] = \deg r - t\sum_{i=1}^d x_i$, where $r$
and $x_i$ are homogeneous elements of $R$.

\begin{defn} 
In the above setting, Goto and Watanabe define the $a$--$invariant$ of $R$ as 
the highest integer $a(R)=a$ such that $[H_m^d(R)]_a$ is nonzero. 
\end{defn}

When $R$ is a ring of characteristic $p$, the Frobenius homomorphism of $R$ 
gives a natural Frobenius action on $H_m^d(R)$ where 
$$
F:[r+(x_1^t,\dots, x_d^t)] \mapsto [r^p+(x_1^{pt},\dots, x_d^{pt})], \ 
\text{ see \cite{FW} or \cite{Sminv}.}
$$
For a graded $R$--module $M$, define 
$M^{(n)}=\oplus_{i \in {\mathbb {Z}}}[M]_{in}$. 
With this notation, it follows from \cite[Theorem 3.1.1]{GW} that 
$$
H_{m_{R^{(n)}}}^d(R^{(n)}) \cong (H_m^d(R))^{(n)}.
$$

The following theorem, \cite[Theorem 7.12]{HHjalg}, indicates the importance of
the $a$--invariant in the study of graded F--rational rings.

\begin{thm} 
A graded \CM normal ring $R$ over a field of prime characteristic $p$ is 
F--rational if and only if $a(R) < 0$ and the ideal generated by some 
homogeneous system of parameters for $R$ is Frobenius closed. 
\label{fratequiv} 
\end{thm}

\section{F--rationality of Veronese subrings}

The following proposition, well-known to the experts, addresses 
the existence of F--rational Veronese subrings. 

\begin{prop}
Let $R$ be a graded \CM domain of dimension $d$, which is
locally F--rational on the punctured spectrum $\Spec R-{m}$. (This is satisfied,
in particular, if $R$ has an isolated singularity.) Then $[H_m^d(R)]_0=0$ if
and only if the Veronese subring $R^{(n)}$ is F--rational for all integers $n
\gg 0$. In particular if $a(R)<0$, then $R^{(n)}$ is F--rational for all
integers $n\gg 0$. 
\end{prop}

\begin{proof}
Note that we have $[H_m^d(R)]_0 \subseteq 0_{H_m^d(R)}^*$, since for $z\in
[H_m^d(R)]_0$ we get $cz^q=0$ for all $q=p^e$, when $c \in m$ is of a
sufficiently large degree. Consequently if $R^{(n)}$ is F--rational for some
$n$, we must have $a(R^{(n)})<0$, but then $[H_m^d(R)]_0=0$.

For the converse first note that since $R$ is F--rational on the punctured
spectrum, Theorem \ref{longlist} (6) says that $0_{H_m^d(R)}^*$ must be killed
by a power of the maximal ideal $m$, and so is of finite length. As
$[H_m^d(R)]_0=0$, for large positive integers $n$ we see that
$H_{m'}^d(R^{(n)}) \cong (H_m^d(R))^{(n)}$ contains no nonzero element of
$0_{H_m^d(R)}^*$ where $m'$ denotes the homogeneous maximal ideal of $R^{(n)}$.
If $u \in 0_{H_{m'}^d(R^{(n)})}^*$ then $u \in 0_{H_m^d(R)}^* \cap
H_{m'}^d(R^{(n)})$ and so $u = 0$. Hence $R^{(n)}$ is F--rational for $n \gg
0$. \end{proof}

\begin{ex} Let $R=K[X,Y,Z]/(X^2+Y^3+Z^5)$ where $K$ is a field of prime
characteristic $p$. We make this a graded ring by setting the weights of $x$,
$y$ and $z$ to be $15$, $10$ and $6$ respectively. We determine the positive
integers $n$ for which the Veronese subring $R^{(n)}$ is F--rational. This
shall, of course, depend on the characteristic $p$ of $R$.

First note that $a(R) =-1$ with this grading. If $p \ge 7$, it is easy to
verify that the ring $R$ is F--regular. Consequently every Veronese subring of
$R$, being a direct summand of $R$, is also F--regular. For $p=2$, $3$ or $5$,
$x^p \in (y^p,z^p)$, and so $R$ is not F--rational. It is also easily checked
that the action of the Frobenius on $H_m^2(R)$ is injective in degree $\le -2$
with the one exception of $p=2$ where elements in degree $-7$ are mapped to
zero under the action of the Frobenius, specifically $F(xy^{-1}z^{-2})=0$ in
$H_m^2(R)$. Recall that $H_{m_{R^{(n)}}}^2(R^{(n)})$ is generated by elements
of $H_m^2(R)$ whose degree is a multiple of $n$. Consequently for $n \ge 2$ the
action of the Frobenius on $H_{m_{R^{(n)}}}^2(R^{(n)})$ is injective, with the
one exception. Using the arguments in the proof of the above proposition, we see
that $R^{(n)}$ is F--rational for all $n \ge 2$, excluding the case when $p=2$
and $n=7$. 
\label{fratex} 
\end{ex}

\section{Rational coefficient Weil divisors}

We review some notation and results from \cite{De}, \cite{Wadem} and
\cite{Wadim2}, as well as make a few observations which we shall find useful
later in our study.

\begin{defn} 
By a {\it rational coefficient Weil divisor}\/ (or a ${\mathbb Q}$--{\it
divisor}) on a normal projective variety $X$, we mean a ${\mathbb Q}$--linear 
combination of codimension one irreducible subvarieties of $X$. For $D = \sum
n_iV_i$ with $n_i \in {\mathbb Q}$, we set $[D]= \sum [n_i]V_i$, where $[n]$
denotes the greatest integer less than or equal to $n$, and define
${\mathcal{O}}_X(D) ={\mathcal{O}}_X([D])$.

Let $D=\sum(p_i/q_i)V_i$ where the integers $p_i$ and $q_i$ are relatively
prime and $q_i > 0$. We define $D' = \sum((q_i-1)/q_i)V_i$ to be the {\it
fractional part\/} of $D$. Note that with this definition of $D'$ we have
$-[-nD] = [nD+D']$ for any integer $n$. 
\end{defn}

Given an ample ${\mathbb Q}$--divisor $D$ (i.e., such that $ND$ is an ample 
Cartier divisor for some $N \in {\mathbb N}$), we construct the 
{\it generalized section ring}:
$$
R = R(X,D) = \oplus _{n \ge 0} H^0(X,{\mathcal {O}}_X(nD))T^n \subseteq
K(X)[T].
$$
With this notation, Demazure's result (\cite[3.5]{De}) is:

\begin{thm} 
Let $R=\oplus _{n \ge 0}R_n$ be a graded normal ring. Then there
exists an ample ${\mathbb Q}$--divisor $D$ on $X = \Proj R$ such that
$$
R = \oplus _{n \ge 0} H^0(X,{\mathcal {O}}_X(nD))T^n \subseteq K(X)[T],
$$
where $T$ is a homogeneous element of degree one in the quotient field of $R$.
\label{demazure}
\end{thm}

\begin{ex}
Take the ${\mathbb Q}$--divisor 
$$
D=(-1/2)V(S)+(1/3)V(T)+(1/5)V(S+T)
$$ on ${\mathbb {P}}^1 = \Proj K[S,T]$ where $V(S)$, e.g., denotes the 
irreducible subvariety defined by the vanishing 
of $S$. Fix $T$ as the degree one element. Then 
$$
R = \oplus_{n \ge 0} H^0({\mathbb {P}}^1,{\mathcal {O}}_{{\mathbb {P}}^1}(nD))T^n 
 = K[X,Y,Z]/(X^2+Y^3+Z^5), \ \text{ where}
$$
$X=(S^8T^{10})/(S+T)^3, Y=(S^5T^7)/(S+T)^2$, and $Z=(-S^3T^4)/(S+T)$. 
\label{favouriteex}
\end{ex}

\begin{rem} 
Let $R = R(X,D)$ be as above. Then the Veronese subring $R^{(n)}$ is given by
$R^{(n)} \cong R(X,nD)$. For a rational function $f \in K(X)$ we have an
isomorphism $R(X,D) \cong R(X,\div(f)+D)$. If $R$ is generated over $K$ by its
elements of degree one, we have $R = R(X,[D])$. Note that $[D]$ is a Weil
divisor, i.e., has integer coefficients. 
\label{qdiv} 
\end{rem}

\section{Results in dimension two}

In the following theorem, we summarize some familiar results about
graded rings of dimension two. 

\begin{thm} 
Let $R$ be a graded normal ring of dimension two, which is
generated by degree one elements over an algebraically closed field.
Then the following statements are equivalent:

\item $(1)$\quad $R$ is isomorphic to a Veronese subring of a polynomial
ring in two variables.

\item $(2)$\quad $R$ is F--regular.

\item $(3)$\quad $R$ is F--rational.

\item $(4)$\quad$R$ has a negative a--invariant.
\label{2dim1form}
\end{thm}

\begin{proof}
The implications (1) $\Rightarrow (2) \Rightarrow (3) \Rightarrow (4)$
follow easily. For $(4) \Rightarrow (1)$ note that $X=\Proj R$ is a
nonsingular projective curve. Since $[H_m^2(R)]_0 = 0$, we have
$H^1(X,{\mathcal {O}}_X)=0$ and so $X$ is of genus zero, i.e., ${\mathbb P}^1$.
Consequently $R \cong R({\mathbb P}^1,D)$ where $D$ is a Weil divisor on
${\mathbb P}^1$. Hence $D$ is linearly equivalent to ${\mathcal O}(m)$ for
some $m \in {\mathbb N}$ and $R \cong R({\mathbb P}^1, {\mathcal O}(m) )
\cong (K[X_0,X_1])^{(m)}$. \end{proof}

As an easy consequence of the above, we have:

\begin{thm}
Let $R$ be a graded domain of dimension two, with an
isolated singularity, which is finitely generated over an algebraically
closed field. If $a(R) < 0$, there exists a positive integer $n$ such
that $R^{(n)}$ is isomorphic to a Veronese subring of a polynomial ring in
two variables over $K$. In particular, some Veronese subring of $R$ is
F--regular. 
\label{2dimver} 
\end{thm}

\begin{proof}
Note that $R$ is excellent and so $R'$, the integral closure $R$ in its
fraction field, is module--finite over $R$. Since $R$ has an isolated
singularity, the conductor (i.e., the largest common ideal of $R$ and $R'$)
is primary to the maximal ideal of $R'$, by which $R_i = R_i'$ for all $i
\gg 0$. We may therefore choose a positive integer $k$ such that $R^{(k)}$
is normal, and then choose an appropriate multiple $n$ of $k$, by Lemma
\ref{existver}, such that $R^{(n)}$ is generated by elements of equal
degree. We are now in a position to apply the above theorem to conclude
that $R^{(n)}$ is isomorphic to a Veronese subring of a polynomial ring in
two variables. 
\end{proof}

\begin{ex}
Let $S=K[X, \ Y, \ Z]/(X^3-YZ(Y+Z))$ where $K$ is a field of characteristic 
$p \equiv 1 \pmod 3$ and consider the subring
$$
R=K[X, \ Y^3, \ Y^2Z, \ YZ^2, \ Z^3]/(X^3-YZ(Y+Z)).
$$ 
It is proved in \cite{HHjalg} that $R$ is F--rational but 
not F--regular, see also \cite{Wadim2}. Since $R^{(3)}$ is generated by 
elements of equal degree, it must be isomorphic to a Veronese subring of a 
polynomial ring by Theorem \ref{2dim1form}. Indeed,
$$
R^{(3)}= K[Y^3,\ Y^2Z,\ YZ^2,\ Z^3].
$$

\end{ex}

\begin{ex}
Let $R=K[t,t^4x,t^4x^{-1},t^4(x+1)^{-1}]$ where $K$ is a field of prime
characteristic $p$. This is one of the examples in \cite{Wafpure} of 
rings which are F--rational but not F--pure; for a different proof see 
\cite{HHjalg}. By mapping a polynomial ring onto it, we may write $R$ as 
$$
R=K[T,\ U,\ V,\ W]/(T^8-UV, \ T^4(V-W)-VW, \ U(V-W)-T^4W).
$$
This is graded by setting the weights of $t$, $u$, $v$ and $w$ to be $1$, 
$4$, $4$ and $4$ respectively. Note that
$$
R^{(4)}=K[S,\ U,\ V,\ W]/(S^2-UV, \ S(V-W)-VW, \ U(V-W)-SW)
$$
where we relabel $T^4$ as $S$. Then $R^{(4)}$ is generated by elements of 
equal degree, and is isomorphic to $K[X^3,X^2Y,XY^2,Y^3]$ by setting 
$S=XY(X-Y)$, $U=XY^2$, $V=X(X-Y)^2$, and $W=Y(X-Y)^2$. 
\end{ex}

By Theorem \ref{2dimver} we know that a graded normal ring $R$ of dimension
two over an algebraically closed field has a Veronese subring $R^{(n)}$ which
is F--regular. We next show that if $R$ is a hypersurface, there exists $n$
such that $R^{(n)}$ is actually an F--regular hypersurface.

\begin{thm} 
Let $R$ be a graded normal hypersurface of dimension two with 
$a(R) <0$. Then there exists a positive integer $n$ such that the Veronese
subring $R^{(n)}$ is an F--regular hypersurface.
\end{thm} 

\begin{proof}
Let $R=K[X,Y,Z]/(f)$ where $x$, $y$ and $z$ have weights $m$, $n$ and $r$ 
respectively. We may assume without any loss of generality that $m$, $n$ and
$r$ have no common factor. If $d = \gcd (m,n)$, then by our assumption $d$  and
$r$ are relatively prime. Therefore $f$ must be a polynomial in $x$, $y$ and
$z^d$. Consequently $R^{(n)}$ is again a hypersurface, and satisfies all the
initial hypotheses, and so we may assume that $R$ satisfies the extra
hypothesis that $m$, $n$ and $r$ are pairwise relatively prime. Assume further
that $m \ge n \ge r$. We consider the two cases: a) $n=1$ and $r=1$,
and b) $m > n > r$. Note that it suffices to show that $R$ is F--rational,
since it is indeed a hypersurface.

We first eliminate the case (\#)\ when $f$ is of the form $XH(Y,Z)+G(Y,Z)$. We 
may take a system of parameters of $R$ of the form $x$, $t$ where $t$ is the
image in $R$ of a polynomial $T \in K[X,Y,Z]$ involving only $Y$ and $Z$. If
$R$ is not F--rational, then since $a(R)<0$, $(x,t)$ cannot be F--pure. Hence
for some $q=p^e$, we have $s^q \in (x^q,t^q)$ while $s \notin (x,t)$. Again,
we may assume that $s$ is the image in $R$ of a polynomial $S \in K[X,Y,Z]$
involving only $Y$ and $Z$. This means that in $K[X,Y,Z]$, we have $S^q \in
(X^q,T^q, XH+G)$ but then $S^q \in (T^q,G^q)$ and so $S \in (T,G)$ in
$K[X,Y,Z]$, giving us the contradiction $s \in (x,t)$.

a)\quad We have $a(R)=\deg f -(m+n+r)<0$, and so $\deg f < m+2$ since $n=r=1$. 
This forces $f$ to be of the form (\#).

b)\quad Since $a(R)=\deg f -(m+n+r)<0$, we have $\deg f < m+n+r < 3m$. 
Hence up to a scalar multiple, $f$ is of the form $XH(Y,Z)+G(Y,Z)$ or 
$X^2+G(Y,Z)$. Note that the first case has already been handled.

Now suppose $f=X^2+G(Y,Z)$. Then $\deg f=2m < m+n+r$ and so $3 < m < n+r$, 
consequently $G$ cannot involve a term of the form $Y^2Z^l$ where $l \ge 2$. If
$G$ has a term $Y^k$, then $2m=kn$ and so $n=1$ or $2$. Since $n > r$, we can
only have $n=2$ and $r=1$, but this too is impossible. Hence $f$ can only be of
the form $f=X^2+aZ^k+bYZ^l+cY^2Z$ where $a$, $b$ and $c$ are scalars. $R$ is
normal, and so $c$ must be non--zero since $l \ge 2$ and $k \ge 2$. It follows
that $2m=2n+r$. If $a$ is non--zero, $2m=rk$ and since $r$ is even, we can only
have $r=2$. But then $m=n+1$, and so $r$ divides either $m$ or $n$, a
contradiction. Hence $a=0$, and so $f=X^2+bYZ^l+cY^2Z$. If $b$ were non--zero,
then we would have $n+rl=2n+r$, i.e., $n=r(l-1)$, which forces $r=1$. However
we know $r$ to be even, and so $b=0$. We are left with $f=X^2+cY^2Z$ but this
is ruled out since $R$ is normal. 
\end{proof}

\section{F--regular Veronese subrings}

We begin by recalling a theorem of Watanabe, \cite[Theorem 3.4]{Wadim2}:

\begin{thm} 
Let $D_1$ and $D_2$ be ample ${\mathbb Q}$--divisors on a normal projective
variety $X$. If the fractional parts $D_1'$ and $D_2'$ are equal, then the ring
$R(X,D_1)$ is F--regular (F--pure) if and only if the ring $R(X,D_2)$ is 
F--regular (F--pure). 
\label{fracpart} 
\end{thm}

A complete proof of the theorem, as stated above, relies on the
characterization of strong F--regularity in terms of the tight closure of the 
zero submodule of the injective hull of the residue field, \cite[Proposition
7.1.2]{Smthesis}, as well as the results of \cite{LS}.

\begin{cor} 
Let $R$ be a graded normal ring which is generated by degree
one elements over a field. Then either $R$ is F--regular (F--pure), or else no 
Veronese subring of $R$ is F--regular (F--pure). 
\label{doesntgetbetter}
\end{cor}

\begin{proof}
Since $R$ is generated by its elements of degree one, we have $R = R(X,D)$, 
where $D$ is a Weil divisor, i.e., has $D'=0$. Also, $(nD)'=0$ where $n$ is any 
positive integer. By the above Theorem, $R=R(X,D)$ is F--regular (F--pure)
if and only if $R^{(n)} \cong R(X,nD)$ is F--regular (F--pure). 
\end{proof}

As an application of this result, we now construct a family of rings with
negative $a$--invariants, which have no F--pure Veronese subrings. This
shows that a result corresponding to Theorem \ref{2dimver} is no longer
true in higher dimensions.

\begin{ex}
Let $R=K[X_0, \dots, X_d]/(X_0^3+ \dots +X_d^3)$ with $d \ge 3$, where $K$ is a
field of characteristic $2$. It is readily seen that 
$x_0^2 \in (x_1, \dots, x_d)^*$, since $x_0^4
\in (x_1, \dots, x_d)^{[2]}$. Hence $R$ is not F--pure, and since it is
generated by elements of degree one, Corollary \ref{doesntgetbetter} shows
that $R$ has no F--regular or F--pure Veronese subrings. Note that $a(R)=2-d<0$.

We can also see that $R^{(n)}$ is not F--pure (for any $n > 0$) by showing that
the element $x_0^d(x_1 \dotsm x_d)^{n-1}$ is in the Frobenius closure of the 
ideal 
$$ 
(x_0^{d-2} x_1^n x_2^{n-1} \dotsm x_{d-1}^{n-1},
\ \ x_0^{d-2} x_2^n x_3^{n-1} \dotsm x_d^{n-1}, 
\ \ \dots, 
\ \ x_0^{d-2} x_d^{n}x_1^{n-1} \dotsm x_{d-2}^{n-1}), 
$$
although not in the ideal itself.

For all $n \ge 2$, the ring $R^{(n)}$ is an example of a graded ring
generated by degree one elements (with an isolated singularity and a
negative $a$--invariant) which is F--rational but not F--pure.
\label{nover} 
\end{ex}
 
\begin{rem} 
The examples above are not completely satisfactory as they are not valid in
the characteristic zero setting: in fact, for $d \ge 3$, the ring 
$R={\mathbb Q}[X_0, \dots, X_d]/(X_0^3+ \dots +X_d^3)$ is of F--regular type.
Characteristic zero examples turn out to be much more subtle, and we construct 
these in the next section. 
\end{rem}

We again return to the ring $R=K[X,Y,Z]/(X^2+Y^3+Z^5)$, and this time
determine its F--regular and F--pure Veronese subrings.

\begin{ex} 
Let $R=K[X,Y,Z]/(X^2+Y^3+Z^5)$ where $K$ is a field of prime characteristic
$p$, and the grading is as before. For $p \ge 7$ the ring $R$ is F--regular,
and therefore so is any Veronese subring $R^{(n)}$. We now determine when
$R^{(n)}$ is F--regular assuming $p$ is either $2$, $3$ or $5$.

Note that the Veronese subrings $R^{(2)}$, $R^{(3)}$ and $R^{(5)}$ are in
fact polynomial rings. Therefore when $n$ is divisible by one of $2$, $3$
or $5$, $R^{(n)}$ is a direct summand of a polynomial ring, and so is
F--regular. We show that these are the only instances when $R^{(n)}$ is
F--regular, or even F--pure.

Recall from Example \ref{favouriteex} that $R = R(X,D)$ where $X=\Proj K[S,T]$
and $D=(-1/2)V(S)+(1/3)V(T)+(1/5)V(S+T)$. If $n$ is relatively prime to $30$,
the ${\mathbb Q}$--divisor $nD$ has the same fractional part as $D$, and so
$R^{(n)} \cong R(X,nD)$ is not F--pure or F--regular by Theorem \ref{fracpart}.
 
We can also construct explicit instances of Frobenius closure to
illustrate why $R^{(n)}$ is not F--pure when $n$ is relatively prime to
$30$. Since $n$ is relatively prime to the weight of $y$, the ring
$R^{(n)}$ has a unique monomial of the form $xy^l$ with $0 < l < n$.
Similarly there is a unique integer $m$ with $0 \le m < n$ such that 
$y^{l+1}z^m \in R^{(n)}$, and a unique integer $r$ with $0 < r < n$ such that 
$x^r z \in R^{(n)}$. We claim that 
\begin{align*} 
x^{r+1}y^{rl+l}z^d  \ &\in \ (x^ry^{rl+l+1}z^m, \ x^rz^{d+1})^F , \ \text{ and}\\
x^{r+1}y^{rl+l}z^d  \ &\notin \ (x^ry^{rl+l+1}z^m, \ x^rz^{d+1}). 
\end{align*} 
The second statement is true in $R$ and so also in $R^{(n)}$, while the 
first assertion follows from
$$ 
(x^{r+1}y^{rl+l}z^d)^p \ \in \ ((x^ry^{rl+l+1}z^m)^p, \ 
(x^rz^{d+1})^p) \ \text{ for }p=2, \ 3 \text{ or }5. 
$$ 
\end{ex}

\begin{ex}
We saw that the F--purity and F--regularity of a ring $R=R(X,D)$ depend only on 
the fractional part $D'$ of the ${\mathbb Q}$--divisor $D$. This is by no means
true of F--rationality and F--injectivity (i.e., the injectivity of the
Frobenius action on the highest local cohomology module). As an example of this,
consider the ${\mathbb Q}$--divisors on $\Proj K[S,T]$
\begin{align*}
E &= (1/2)V(S)+(1/3)V(T)+(1/5)V(S+T) \ \text{ and} \\
D &= (-1/2)V(S)+(1/3)V(T)+(1/5)V(S+T).
\end{align*}
which have the same fractional part. Then
$$
S = \oplus _{n \ge 0} H^0(X,{\mathcal {O}}_X(nE))T^n \cong K[A,B,C,T]/I
$$
where $I = (AB-T^5, \ BC+CT^3-BT^5, \ AC+CT^2-ABT^2)$ and $A=T^3/S$, $B=ST^2$
and $C=ST^5/(S+T)$. If the characteristic of $K$ is $2$, $3$ or $5$, the ring
$R=R(X,D)=K[X,Y,Z]/(X^2+Y^3+Z^5)$ is not F--rational (or F--injective) as we
saw in Example \ref{fratex}. We claim that the ring $S$ is however F--rational.
To see this note that $a(S) < 0$, and so it suffices by Theorem \ref{fratequiv}
to verify that the ideal $I$ generated by the homogeneous system of parameters
$t$, $a^{15} + b^{10} + c^6$ is Frobenius closed. However this is easily
verified: the ring $S/tS \cong K[A,B,C]/(AB,BC,CA)$ is F--pure since the ideal
$(AB,BC,CA)$ is generated by square free monomials, see \cite[Proposition
5.38]{HRcohom}. 
\end{ex}

\begin{rem} 
Let $R$ be a \CM ring with an isolated singularity, which is generated by
degree one elements over an algebraically closed field. For a two dimensional
ring $R$, a negative $a$--invariant forces $R$ to be F--regular, although for
rings of higher dimensions this is no longer true: in Example \ref{nover} we
constructed rings $R$ of dimension $d>3$, with $a(R) = 2-d$, which were not
F--regular. Smith has pointed out that if $R$ satisfies the stronger condition
that $a(R) \le 1-d$, then $\Proj R$ is a variety of minimal degree. These are
completely classified (see, for example, \cite{EHmin}) and it is easily
verified that in this case $R$ is F--regular, see \cite[Remark 4.3.1]{Smgraded}. 
\end{rem}

\section{The case of characteristic zero}

Hochster and Huneke have defined analogous notions of tight closure for 
rings essentially of finite type over a field of characteristic zero, see 
\cite{HHjams, HHchar0}. However we can also define notions corresponding to
F--regularity, F--purity, and F--rationality in characteristic zero, without 
using a closure operation. 

Consider the ring $R=K[X_1, \dots, X_n]/I$ where $K$ is a field of
characteristic zero. Choose a finitely generated $\mathbb Z$--algebra $A$ such
that $R_A = A[X_1, \dots, X_n]/I_A$ is a free $A$--algebra, with $R \cong R_A
\otimes_A K$. Note that the fibers of the homomorphism $A \to R_A$ over maximal
ideals of $A$ are finitely generated algebras over fields of prime
characteristic. 

\begin{defn} 
Let $R$ be a ring finitely generated over a field of characteristic
zero. Then $R$ is said to be of {\it F--regular type}\/ if there exists a
finitely generated $\mathbb Z$--algebra $A \subseteq K$ and a finitely
generated $A$--algebra $R_A$ such that $R \cong R_A \otimes_A K$, and for all
maximal ideals $\mu$ in a Zariski dense subset of $\Spec A$, the fiber rings
$R_A \otimes_A A/\mu$ are F--regular.

Similarly, $R$ is said to be of {\it F--pure type}\/ if for all maximal ideals 
$\mu$ in a Zariski dense subset of $\Spec A$, the fiber rings 
$R_A \otimes_A A/\mu$ are F--pure.
\end{defn}

\begin{rem}
Some authors use the term F--pure type (F--regular type) to mean that 
$R_A \otimes_A A/\mu$ is F--pure (F--regular) for all maximal ideals $\mu$ in a 
Zariski dense {\it open}\/ subset of $\Spec A$. 
\end{rem}

All our positive results towards the existence of F--rational and F--regular
Veronese subrings in prime characteristic do have corresponding statements in
the characteristic zero situation. However we have so far not exhibited a
normal Cohen--Macaulay ring, generated by degree one elements over a field of
characteristic zero, which has an isolated singularity and a negative
$a$--invariant but is not of F--regular type. N.~Hara has pointed out to us a
geometric argument for the existence of such rings using a blow--up of
${\mathbb P}^2$ at nine points. In this section, we construct a large family of
explicit examples of such rings of dimension $d \ge 3$.

\begin{ex}
Take two relatively prime homogeneous polynomials $F$ and $G$ of degree $d$ in
the ring ${\mathbb Z}[X_1, \dots, X_k]$, where $k \ge 3$, such that $G$ is
monic in $X_k$ and the monomial $X_k^d$ does not occur in $F$. Using $F$ and
$G$, construct the hypersurface $S={\mathbb Q}[S,T,X_1, \dots, X_k]/(SF-TG)$
and let $R$ be the subring of $S$ generated by the elements $sx_1, \ \dots, \
sx_k, \ tx_1, \ \dots, \ tx_k$.

For suitably general choices of the polynomials $F$ and $G$ of degree $d=k$ the
ring $R$ has only isolated singularities, and we show that it is \CM with $a(R)
= -1$, and is not of F--regular type. For an explicit example, take $k=3$,
$F=X_1X_2X_3$ and $G=X_1^3 + X_2^3 + X_3^3$.

We shall prove that $R$ is \CM whenever $d \le k$. We first show that the
Hilbert polynomial multiplicity of $R$ is $d(k-1)+1$, and then construct a
system of parameters such that the ring obtained by killing this system of
parameters has length $d(k-1)+1$.

We construct a basis for the vector space generated by the monomials of degree
$n \gg 0$, $s^it^{n-i}x_1^{j_1}x_2^{j_2} \dotsm x_k^{j_k}$, where the $j_r$ are
nonnegative integers which add up to $n$. The relations permit us to express 
$tx_k^d$ in terms of other monomials. Let $[u_1, \dots, u_m]^i$ denote the set 
$\mathcal S$ of monomials of degree $i$ in $u_1, \dots, u_m$, and for two such sets, 
let $\mathcal S \cdot T$ denote the product of all possible pairs from $\mathcal S$ and 
$\mathcal T$. In this notation, for $n \gg 0$, the following monomials constitute a 
basis for $R_n$: 
\begin{align*} 
& [s,t]^n \cdot [x_1, \dots, x_{k-1}]^n, \\
& [s,t]^n \cdot [x_1, \dots, x_{k-1}]^{n-1} \cdot [x_k], \\
& \hdots \\
& [s,t]^n \cdot [x_1, \dots, x_{k-1}]^{n-d+1} \cdot [x_k]^{d-1}, \\
& [s]^n \cdot [x_1, \dots, x_k]^{n-d} \cdot [x_k]^d.
\end{align*}
Consequently for large $n$ the vector space dimension of $R_n$ is
$$ 
(n+1)\left\{ \binom{n+k-2}{k-2} + \dots + 
\binom{n-d+1+k-2}{k-2} \right \}+ \binom{n-d+k-1}{k-1}. 
$$ 
As a polynomial in $n$, the leading term of this expression is 
$$ 
n\left\{ \frac{n^{k-2}}{(k-2)!} + \dots + \frac{n^{k-2}}{(k-2)!} \right\} + 
\frac{n^{k-1}}{(k-1)!} = \frac{n^{k-1}(d(k-1)+1)}{(k-1)!}, 
$$ 
and so the Hilbert polynomial multiplicity of $R$ is $d(k-1)+1$. 

The sequence of elements $sx_1, \ sx_2-tx_1, \ sx_3-tx_2, \ \dots, \
sx_k-tx_{k-1}$ is a system of parameters for $R$. Since we have already
verified that the Hilbert polynomial multiplicity of $R$ is $d(k-1)+1$, to
prove that $R$ is Cohen--Macaulay when $d \le k$, it suffices to show that the
length of the ring $T$ obtained by killing this system of parameters is at most
$d(k-1)+1$.

Relabel the generators of $T$ as 
$a_2 = sx_2, \ a_3=sx_3, \ \dots, \ a_k=sx_k, \ a_{k+1}=tx_k$. 
Note that the relations amongst the $a_i$ include the size two minors of the 
matrix 
$$
\begin{pmatrix}
0 & a_2 & \hdots & a_{k-1} & a_k \\
a_2 & a_3 & \hdots & a_k & a_{k+1} 
\end{pmatrix}.
$$
Consequently a generating set for $[T]_{<d}$ is given by
\begin{alignat*} 2 
 \text{deg } & 0 & &: 1, \\
 \text{deg } & 1 & &: a_2, \ a_3, \ \dots, \ a_{k+1}, \\
 \text{deg } & 2 & &: a_2a_{k+1}, \ a_3a_{k+1}, \ \dots, \ a_{k+1}^2, \\
 \text{deg } & 3 & &: a_2a_{k+1}^2, \ a_3a_{k+1}^2, \ \dots, \ a_{k+1}^3, \\
 & & & \dots \\
 \text{deg } & d-1 & &: a_2a_{k+1}^{d-2}, \ a_3a_{k+1}^{d-2},\ \dots, \ a_{k+1}^{d-1}.
\end{alignat*}
In degree $d$ the ring $T$ has $d$ additional independent relations coming
from the equations $s^it^{d-i}f - s^{i-1}t^{d-i+1}g$, for $1 \le i \le d$.
Consequently we need $k-d$ generators for the degree $d$ piece of $T$, and one 
can check that there are no nonzero elements in degree $d+1$. Hence the length 
of $T$ is bounded by $d(k-1)+1$, and this completes the proof that $R$ is
Cohen--Macaulay.

It only remains to show that $R$ is not of F--regular type when $k\le d$.
Consider the fiber $A$ of the map ${\mathbb Z} \to R_{{\mathbb Z}}$ over an
arbitrary closed point $p{\mathbb Z}$. Then $A$ is a finitely generated algebra
over the finite field ${\mathbb Z}/p{\mathbb Z}$, and it suffices to show that 
$A$ is not F--regular. Take the ideal 
$$
I=(sx_1, sx_2, \dots, sx_{k-1}, tx_1, tx_2, \dots, tx_{k-1} )A.
$$ 
It is easily verified that $(tx_k)^{d-1} \notin I$, and we show that
$(tx_k)^{d-1} \in I^*$. 

To see $(tx_k)^{d-1} \in I^*$ it suffices to check that $\alpha_q =
(tx_k)^{(d-1)(q+1)} \in I^{[q]}$. Using the relation $t^d g - t^{d-1}sf$ where 
$1 \le i \le d$, we may rewrite $\alpha_q$ with lower powers of $x_{k}$
occurring in the expressions involved. We can proceed in this manner till we
are left with terms which involve powers of $x_k$ not greater then $d-1$. Hence
$\alpha_q$ is a sum of terms which are multiples of

$$
s^i t^{q(d-1)-i} x_1^{j_1}x_2^{j_2} \dotsm x_{k-1}^{j_{k-1}}, 
\text{ \ where } i \le q(d-1), \text { \ and } \sum_{r=1}^{k-1} j_r = q(d-1).
$$
If $\alpha_q \notin I^{[q]}$, then $j_r < q$ for $1 \le r \le k-1$. 
However on summing these inequalities we get $q(d-1) < q(k-1)$, a
contradiction. 

\end{ex}

\begin{rem}
Consider the polynomial ring $K[X_1, \dots, X_k]$ where $k \ge 3$. It is worth 
noting that the ring $R$, as above, is isomorphic to a subring of 
$K[X_1, \dots, X_k]$, 
$$
R=K[X_1F, \ X_2F, \ \dots, \ X_kF, \ X_1G, \ X_2G, \ \dots, \ X_kG].
$$
We can show that $R$ is \CM precisely when the degree $d$ of $F$
and $G$ is less than or equal to $k$. It would certainly be interesting to 
explore generalizations of this construction. 
\end{rem}

\section*{Acknowledgments}

The author wishes to thank Melvin Hochster and Karen Smith for many valuable
discussions.

\end{document}